\documentclass[11pt,a4paper]{amsart}
\usepackage{amsmath,amsfonts,amsthm,amssymb,amscd}
\usepackage{appendix}
\usepackage[latin1]{inputenc}
\usepackage{mathrsfs}
\usepackage{latexsym}
\usepackage{graphicx}
\usepackage{enumerate}
\def\Z{\mathbb{Z}}
\def\M{\mathbb{Z}_{\ge{0}}}
\def\N{\mathbb{Z}_{>0}}
\def\R{\mathbb{R}}
\def\W{\mbox{WN}}

\newcommand{\mex}{\operatorname{mex}}
\newcommand\T{\rule{0pt}{2.6ex}}

\theoremstyle{definition}
\newtheorem{Thm}{Theorem}[section]
\newtheorem{Prop}[Thm]{Proposition}
\newtheorem{Lemma}[Thm]{Lemma}

\newtheorem{Def}{Definition}

\newtheorem{Ex}{Example}
\newtheorem{Rem}{Remark}

\begin{document}
\title{Restrictions of $m$-Wythoff Nim and $p$-complementary Beatty Sequences}

\author{Urban Larsson}
\email{urban.larsson@chalmers.se}
\address{Mathematical Sciences, Chalmers University of Technology
and University of Gothenburg, G\"oteborg, Sweden}
\date{\today }

\begin{abstract}
Fix a positive integer $m$. The game
of \emph{$m$-Wythoff Nim} (A.S. Fraenkel, 1982) is
a well-known extension of \emph{Wythoff Nim}, a.k.a 'Corner the Queen'.
Its set of $P$-positions may be represented by a pair of increasing
sequences of non-negative integers.
It is well-known that these sequences are so-called
\emph{complementary homogeneous}
 \emph{Beatty sequences}, that is they satisfy Beatty's theorem.
For a positive integer $p$, we generalize the solution of $m$-Wythoff
Nim to a pair of \emph{$p$-complementary}---each positive integer
occurs exactly $p$ times---homogeneous Beatty sequences
$a = (a_n)_{n\in \M}$ and $b = (b_n)_{n\in \M}$, which, for all $n$,
satisfies $b_n - a_n = mn$. By the latter property, we show that
$a$ and $b$ are unique among \emph{all} pairs of non-decreasing
$p$-complementary sequences. We prove that such pairs can 
be partitioned into $p$ pairs of complementary Beatty sequences.
Our main results are that 
$\{\{a_n,b_n\}\mid n\in \M\}$ represents the solution to three
new '$p$-restrictions' of $m$-Wythoff Nim---of
which one has a \emph{blocking maneuver} on the \emph{rook-type} options.
C. Kimberling has shown that the solution of Wythoff Nim satisfies
the \emph{complementary equation} $x_{x_n}~=~y_n~-~1$. We generalize
this formula to a certain '$p$-complementary equation' satisfied by our
pair $a$ and $b$. We also show that one may obtain our new pair of
sequences by three so-called \emph{Minimal EXclusive} algorithms. We 
conclude with an Appendix by Aviezri Fraenkel.

\end{abstract}
\maketitle
\vskip 30pt
\section{Introduction and notation}
The combinatorial game of
\emph{Wythoff Nim} (\cite{Wyt07}) is a so-called (2-player) impartial
game played on two piles of tokens. (For an
introduction to impartial games see \cite{AlNoWo07, BeCoGu82, Con76}.)
As an addition to the rules of the
game of Nim (\cite{Bou02}), where the players alternate in removing any finite
number of tokens from precisely one of the piles (at most the whole pile),
Wythoff Nim also allows removal of the same number of tokens
from both piles. The player who removes the last token wins.

This game is more known as 'Corner
the Queen', invented by R. P. Isaacs (1960), because the game can
be played on a (large) Chess board with one single Queen.
Two players move the Queen alternately but with the restriction
that, for each move, the ($L^1$) distance to the lower left corner,
position $(0, 0)$, must decrease. (The Queen must at all times remain
on the board.) The player who moves to this \emph{final/terminal}
position wins.

In this paper we follow the convention to denote our
players with the \emph{next player} (the player who is in turn to
move) and the \emph{previous player}.
A \emph{$P$-position} is a position from which the previous player can
win (given perfect play).
An \emph{$N$-position} is a position from which the next player can
win. Any position is either a $P$-position or an $N$-position. We
denote \emph{the solution}, the set of
all $P$-positions, of an impartial game $G$,
by $\mathcal{P} = \mathcal{P}(G)$
and the set of all $N$-positions by $\mathcal{N}=\mathcal {N}(G)$.
The positive integers are denoted by $\N$ and the non-negative integers
by $\M$. Let $x = (x_i)$ denote an integer sequence over some index set 
and let $\xi \in \Z$.
Then, we let $x_{ > \xi}$ denote $(x_i)_{i > \xi}$.

\subsection{Restrictions of $m$-Wythoff Nim}
Let $m\in \N$. We next turn to a certain $m$-extension of Wythoff Nim,
studied in \cite{Fra82} by A.S. Fraenkel.
In the game of \emph{$m$-Wythoff Nim}, or just $m\W$ (our notation),
the Queen's 'bishop-type' options are extended so
that $(x + i,y + j)\rightarrow(x,y)$ is legal if $\mid\! i-j\!\mid\; < m$,
$i, j, x, y\in \M $, $i > 0$ or $j > 0$.
The \emph{rook-type} options are as in Nim. Hence $1$-Wythoff Nim is
identical to Wythoff Nim.

In this paper we define three new \emph{restrictions}
of $m$-Wythoff Nim---here a rough outline:

\begin{itemize}
\item The first has a so-called
\emph{blocking maneuver}/\emph{Muller Twist}
(see also \cite{HoRe, SmSt02} and Section 1.2 of this paper) on the
rook-type options, before the next player moves, the previous player
may announce some of these options as forbidden;
\item The second has a certain congruence restriction on the rook-type
options;
\item For the third, a rectangular piece is removed from the lower left
corner of the game board (including position $(0,0)$).
\end{itemize}
Depending on the particular setup, in addition to $(0,0)$,
the first two games may
have a finite number of final positions of the form
$\{0,x\}$\footnote{For integers $0 \le x \le y$ we use the notation
$\{x, y\}$ whenever $(x,y)$ and $(y,x)$ are considered the same.}, $x\in \N$.
In the third game, there are precisely two final positions $\ne (0,0)$.

\subsection{Beatty sequences and $p$-complementarity}
A (general) \emph{Beatty sequence} is a non-decreasing integer sequence
of the form
\begin{align}\label{genbeat}
(\lfloor n\alpha + \gamma \rfloor),
\end{align}
usually indexed over $\Z,\M $ or $\N$,
where $\alpha$ is a positive irrational and $\gamma\in \R$.
S. Beatty \cite{Bea26} is probably most known
for a (re)\footnote{This theorem
was in fact discovered by J. W. Rayleigh, see \cite{Ray94, Bry03}.}discovery
of (the statement of) the following theorem:
If $\alpha $ and $ \beta $ are positive real numbers such that
$\frac{1}{\alpha } + \frac{1}{\beta} = 1$ then $(\lfloor n\alpha \rfloor)_{\N}$
and $(\lfloor n\beta \rfloor )_{\N}$ partition $\N$ if and only
if they are Beatty sequences. This was proved in \cite{HyOs27}
(see also \cite{Fra82}). A pair of sequences that partition $\N$ ($\M$, $\Z$) 
is usually called \emph{complementary} (see \cite{Fra69, Fra73, Kim07, Kim08}).
Let us generalize this notion.

\begin{Def}\label{def1}
Let $p\in \N$ and $Q,R,S\subset \Z$.
Two sequences $(x_i)_{i\in Q}$ and $(y_i)_{i\in R}$ of non-negative integers
are \emph{$p$-complementary} (on $S$), if, for each $n\in S$,
$$\# \{i\in Q\mid x_i = n\} + \# \{i\in R \mid y_i = n\} = p.$$
\end{Def}

\begin{Rem}
Since the main topics in this paper are three games mostly played on
$\M\times \M$, we often 
find it convenient to use $S = \M$ in Definition \ref{def1}.
Also, for our purposes, it will be convenient use $\M$ or $\N$ as the index 
sets $R$ and $Q$. In the Appendix Aviezri Fraenkel provides an 
alternative approach.
\end{Rem}

We study the Beatty sequences $a = (a_n)_{n\in \M}$ and
$b = (b_n)_{n\in \M}$, where for all $n\in \M$,
\begin{align}\label{A}
a_n &= a_n^{m,p} = \left\lfloor \frac{n\phi_{(mp)}}{p}\right\rfloor
\intertext{ and }
b_n&=b_n^{m,p}= \left\lfloor \frac{n(\phi_{(mp)}+mp)}{p}\right\rfloor,\label{B}
\end{align}
and where
\begin{align}\label{Phix}
\phi_{k}=\frac{2-k +\sqrt{k^2+4}}{2}.
\end{align}
We show that $a$ and $b_{>0}$ are $p$-complementary.

In \cite{Wyt07} W.A. Wythoff proved that the solution of Wythoff Nim
is given by $\{\{a_n^{1,1}, b_n^{1,1} \}\mid n\in \M \}$.
Then in \cite{Fra82} it was shown that the solution of $m$-Wythoff Nim is
$$\{\{a^{m,1}_n, b^{m,1}_n\}\mid n\in \M \}.$$

\subsection{Recurrence}
Let $X$ be a strict subset of the non-negative integers. Then
the \emph{Minimal EXclusive} of $X$ is defined
as usual (see \cite{Con76}): $$\mex X:= \min (\M\!\setminus\! X).$$
For  $n\in \M$ put
\begin{align}
x_n = \mex\{x_i,y_i\mid i\in \{0,1,\ldots , n-1] \} \text{ and }
y_n = x_n + mn.\label{framex}
\end{align}
With notation as in (\ref{framex}), it was proved
in \cite{Fra82} that $(x_n)=(a^{m,1}_n)$ and $(y_n)=(b^{m,1}_n)$.
The minimal exclusive algorithm in (\ref{framex}) gives an exponential
time solution to $m\W$ whereas the Beatty-pair in (\ref{A}) and
(\ref{B}) give a polynomial time ditto.
(For interesting discussions on complexity issues
for combinatorial games, see for example \cite{Fra04, FrPe09}.)
We show that one may, for general $m$ and $p$, obtain $a$ and $b$ by
three minimal exclusive algorithms, which
in various ways generalize (\ref{framex}).

It is well-known that the solution of Wythoff Nim
satisfies the \emph{complementary equation}
(see for example \cite{Kim95, Kim07, Kim08})
$$x_{x_n}~=~y_n~-~1.$$  For arbitrary positive integers
$m$ and $p$, we generalize this formula to a '$p$-complementary equation'
\begin{align}\label{pcomp}
x_{\varphi_n} = y_n - 1,
\end{align}
where $\varphi_n := \frac{x_n + (mp - 1)y_n}{m}$ (or $\varphi _n := py_n - n$), 
and show that a solution is given by $x = a$ and $y = b$. 


\subsection{I.G. Connell's restriction of Wythoff Nim}
In the literature there is another generalization of Wythoff Nim that
is of special interest to us.
Let $p\in \N$. In \cite{Con59} I.G. Connell studies the
restriction of Wythoff Nim, where the the rook-type options are restricted
to jumps of precise multiples of $p$.
This game we call Wythoff modulo-$p$ Nim and denote with $\W^{(p)}$. Hence
Wythoff modulo-1 Nim equals Wythoff Nim.

From \cite{Con59} one may derive
that $\mathcal{P}(\W^{(p)})=\{\{a_n^{1,p}, b_n^{1,p}\}\mid n\in \M\}$.

\begin{table}[h!]
\begin{center}
  \begin{tabular}
{| l || c | c | c | c | c | c | c | c | c | c | c | c | c | c | c | c | c|}
    \hline
    $b^{1,3}_n$\T &0&1&2&4&5&7&8&10&11&12&14&15&17&18&20&21&22 \\
    $a^{1,3}_n$\T &0&0&0&1&1&2&2&3&3&3&4 &4 &5 &5 &6 &6&6\\ \hline
    $n$ \T &0 &1&2&3&4&5&6&7&8&9&10&11&12&13&14&15&16\\
    \hline
  \end{tabular}
\end{center}\caption{The initial $P$-positions of Connell's restriction
of Wythoff Nim, $\W^{(3)}$.}
\end{table}

\begin{figure}[h]\label{Fig1}
\centering
\includegraphics[width=1\textwidth]{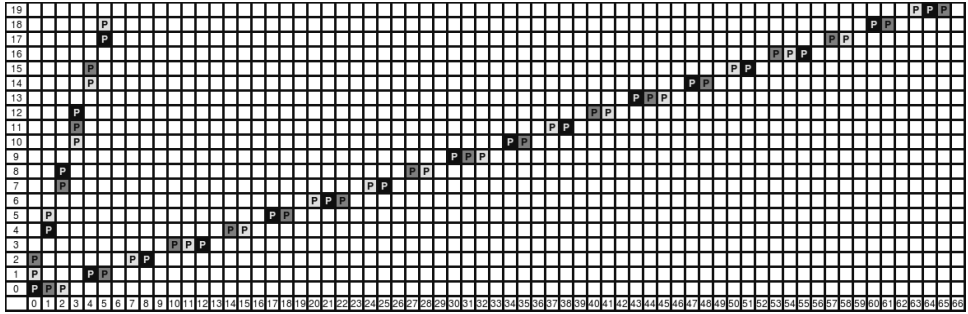}
\caption{The $P$-positions of $\W^{(3)}$
are the positions nearest the origin such that there are precisely
three positions in each row and column and one
position in each NE-SW-diagonal. The black positions represent the first
few $P$-positions of $3$-Wythoff Nim, namely
the positions nearest the origin such
that there is precisely one position in each row and one position in
every third NE-SW diagonal. Positions with lighter shades represent
the solutions of games in our third variation.
}
\end{figure}

\begin{Rem} In Connell's presentation, for
the proof of the above formulas, he uses $p$ pairs of
complementary sequences of integers (in analogy with the discovery of a
new formulation of Beatty's theorem in \cite{Sko57}).
We have indicated this pattern of $P$-positions
with different shades in Figure \ref{Fig1}. In fact,
the black positions, starting with $(0,0)$
are $P$-positions of $3$-Wythoff Nim. More generally, it is immediate by
(\ref{A}) and (\ref{B}) that, for all $p$ and $n$,
$a_n^{p,1} = a^{1,p}_{pn}$ and $b_n^{p,1} = b^{1,p}_{pn}$.
\end{Rem}

\begin{Rem} In \cite{BoFr73}, Fraenkel and I. Borosh study
yet another variation of both $m$-Wythoff Nim and Wythoff modulo-$p$ Nim
which includes a (different from ours) Beatty-type characterization
of the $P$-positions.
\end{Rem}

\subsection{Exposition}
In Section 2, given fixed game constants $m,p\in \N$, we define 
our games, exemplify them and
 state our Main Theorem. Roughly: For each of our games, 
a position is $P$ if and only if it is of the 
form $\{a_n^{m,p}, b_n^{m,p}\}$, with $a$ and $b$ as
in (\ref{A}) and (\ref{B}) respectively.
In Section 3 we generalize Beatty's theorem to pairs of
$p$-complementary sequences, establish that such sequences can 
be partitioned into $p$ pairs of complementary Beatty sequences
and at last prove some arithmetic properties of $a$ and $b$---most important 
of which is that $a$ and $b_{>0}$ are $p$-complementary and, for all $n\in \M$,
satisfy $b_n = a_n + mn$.
Then, in Section 4, we prove that the latter two properties 
make $a$ and $b$ unique among all pairs of non-decreasing sequences.
Section 5 is devoted to equation (\ref{pcomp}) and three
minimal exclusive algorithms.
In Section 6 we prove our game theory results (stated in Section 2) and
 in Section 7 a few questions are posed. At last there is an Appendix, 
provided by Aviezri Fraenkel.

Let us, before we move on to our games, give some more background
to the so-called \emph{blocking maneuver} in the context of Wythoff Nim.

\subsection{A bishop-type blocking variation of $m$-Wythoff Nim}
Let $m, p\in \N$.
In \cite{HeLa06} we gave an exponential time solution
to a variation of $m$-Wythoff Nim with a 'bishop-type' blocking maneuver,
denoted by $p$-Blocking $m$-Wythoff Nim (and with $(m,p)$-Wythoff Nim
in \cite{Lar09}).

The rules are as in $m$-Wythoff Nim, except
that before the next player moves, the previous player is allowed to
block off (at most) $p-1$ bishop-type---note, not $m$-bishop-type---options
and declare that the next player must refrain from these options.
When the next player has moved, any blocking maneuver is forgotten.

The solution of this game is in a certain sense 'very close' to
pairs of Beatty sequences (see also the Appendix of \cite{Lar09}) of the form
\begin{align*}
\left( \left\lfloor n\frac{\sqrt{m^2+4p^2}+2p-m}{2p}\right\rfloor \right)
\text{and}
\left( \left\lfloor n\frac{\sqrt{m^2+4p^2}+2p+m}{2p}\right\rfloor \right).
\end{align*}
But we explain why there can be no Beatty-type solution to
this game for $p>1$.
However, in \cite{Lar09}, for the cases $p\!\mid\!m$, we give a certain
'Beatty-type' characterization. For these kind of questions, see
also \cite{BoFr84}. However, a recent discovery, in \cite{Had, FrPe09},
provides a polynomial time algorithm for the solution of $(m,p)$-Wythoff Nim
(for any combination of $m$ and $p$).

An interesting connection to $4$-Blocking $2$-Wythoff Nim
is presented in \cite{DuGr08}, where
the authors give an explicit bijection of solutions to
a variation of Wythoff's original game,
where a player's bishop-type move is restricted to jumps by
multiples of a predetermined positive integer.

For another variation, \cite{Lar09}
defines the rules of a so-called move-size dynamic variation of two-pile
Nim, $(m,p)$-Imitation Nim, for which the $P$-positions, treated as starting
positions, are identical to the $P$-positions of $(m,p)$-Wythoff Nim.

This discovery of a 'dual' game to $(m, p)$-Wythoff Nim has in its turn
motivated the study of 
'dual' constructions of the rook-type blocking maneuver in this paper.

\section{Three games}
This section is devoted to defining and exemplifying
our game rules and stating our main results. We begin by
introducing some (non-standard) notation whereby we decompose
the Queen's moves into \emph{rook-type} and \emph{bishop-type} ditto.

\begin{Def}\label{moves} 
Fix $m,p \in \N$ and an $l\in \{0, 1, \ldots , m\}$.
\begin{enumerate}[(i)]
  \item An \emph{$(l, p)$-rook} moves as in Nim, but the length
    of a move must be $ip + j > 0$ positions for some $i\in \M$ and
    $j\in \{0,1,\ldots , l-1\}$ (we denote
    a $(0, p)$-rook by a $p$-rook and a $(p, p)$-rook simply by a \emph{rook});
  \item A \emph{$m$-bishop} may move $0 \le i < m$
    rook-type positions and then any number of, say $j\ge 0$,
    bishop-type positions (a \emph{bishop} moves as in Chess),
    all in one and the same move, provided $i + j > 0$ and
    the $L^1$-distance to $(0,0)$ decreases.
\end{enumerate}
\end{Def}

\subsection{Game definitions}
As is clear from Definition \ref{moves} the rook-type options
intersect the $m$-bishop-type options precisely when $m > 1$. For example,
$(2, 3)\rightarrow (1, 3)$ is both a $2$-bishop-type and a rook-type option. 
We will make use of this fact when defining the blocking maneuver.
Therefore, let us introduce some new terminology.

\begin{Def}\label{roob}
Fix $m\in \N$. A rook-type option, which is not of the form of the
$m$-bishop as in Definition \ref{moves} (ii),
is a \emph{roob-type}\footnote{Think of 'roob' as
'ROOk minus $m$-Bishop', or maybe 'ROOk Blocking'} option.
\end{Def}

Hence, for $m=2$, a \emph{roob} may move $(2, 3)\rightarrow (2, 1)$, but not
$(2, 3)\rightarrow (2, 2)$ (both are rook-type options).
Let us define our games.

\begin{Def}\label{games}
 Fix $m, p\in \N$.
\begin{enumerate}
\item [(i)] The game of \emph{$m$-Wythoff $p$-Blocking Nim}, or $m\W^p$,
is a restriction of $m$-Wythoff Nim with a roob-type blocking maneuver.
The Queen moves as in $m$-Wythoff Nim (that is, as the $m$-bishop or the rook),
but with one exception:
Before the next player moves, the previous player may \emph{block off} (at
most) $p-1$ of the next player's roob-type options.
Any blocked option is unavailable for the next player. As usual,
each blocking maneuver is particular to a specific move; that is, when
the next player has moved, any blocking maneuver is forgotten and
has no further impact on the game. (For $p = 1$ 
this game equals $m$-Wythoff Nim.)

\item [(ii)] Fix an integer $0 \le l < p.$ In the
game of \emph{$m$-Wythoff Modulo-$p$ $l$-Nim}, or $m\W^{(l,p)}$,
 the Queen moves as the $m$-bishop
or the $(l, p)$-rook. For $l = 0$ we denote this game by
\emph{$m$-Wythoff Modulo-$p$ Nim} or $m\W^{(p)}$.
(In case $l = p$ the game is simply $m$-Wythoff Nim.)

\item [(iiia)] Fix an integer $0\le l<p.$ In the game of
\emph{$l$-Shifted $m\!\! \times \!\! p$-Wythoff Nim},
or $m\!\!\times \!\!p\W_l$, the Queen moves as in $(mp)$-Wythoff Nim
(that is, as the $(mp)$-bishop or the rook),
except that, if $l > 0$, it is not allowed to move to a position of the form
$(i,j)$, where $0\le i < ml$ and  $0 \le j < m(p-l)$. Hence, for this case,
the terminal positions are $(ml, 0)$ and $(0, m(p-l))$.\footnote{One might
want to think of the game board as if a rectangle with circumference $2mp$
is cut out from its lower left corner.
By symmetry, there are $\lceil \frac{p}{2}\rceil$ rectangle shapes
but (given a starting position) $p$ distinct game boards, see also (3b).
Of course, if cutting out a corner of the nice game board
does not appeal to the players, one might just as well
define all positions inside the rectangle as $N$. Notice the 
close relation of these games to the mis\`{e}re version of $m$-Wythoff Nim 
studied in \cite{Fra84}.
} On the other hand $m\!\times \!p\W_0$ is identical to $(mp)$-Wythoff Nim.

\item [(iiib)]
The game of $m\!\! \times \!\! p$-Wythoff Nim, $m\!\!\times \!\!p\W$:
Before the first player moves, the second player
may decide the parameter $l$ as in (3a). Once the
parameter $l$ is fixed, it remains
the same until the game has terminated, so that for the remainder
of the game, the rules are as in $m\!\!\times \!\!p\W_l$.
\end{enumerate}
\end{Def}

\subsection{Examples}
Let us illustrate some of our games, where our
players are \emph{Alice} and \emph{Bob}---Alice
makes the first move (and Bob makes the first
blocking maneuver in case the game has a Muller twist).

\begin{Ex}\label{ex1}
Suppose the starting position is $(0,2)$ and the game is $2\W^2$. Then
the only bishop-type move is $(0,2)\rightarrow (0,1)$. There is
precisely one roob-type option, namely $(0,0)$. Since this is a terminal
position Bob will block it off from Alice's options, so that
Alice has to move to $(0,1)$.
The move $(0,1)\rightarrow (0,0)$ cannot be blocked off for
the same reason, so Bob wins. If $y\ge 3$ there is always a
move $(0,y)\rightarrow (0,x)$, where $x=0$ or $2$. This is because the
previous player may block off at most one option. Altogether, this gives that
$\{0,y\}$ is $P$ if and only if $y=0$ or $2$.
\end{Ex}
\begin{Ex}\label{ex2}
Suppose the starting position is $(0,2)$ and the game is $2\W^{(2)}$.
Alice can move to $(0,0)$, since $0\equiv 2 \pmod 2$, so $(0,2)$ is $N$. On
the other hand, the position $(0,3)$ is $P$ since the only options
are $(0,2)$ and $(0,1)$. (The latter is $N$ since the $2$-bishop
can move $(0,1)\rightarrow (0,0)$.)
\end{Ex}
\begin{Ex}\label{ex3}
Suppose the starting position is $(0,2)$ and the game is $2\W^{(2,4)}$.
Alice cannot move to a $P$-position since $2 - 0\not \equiv 0,1 \pmod 4$ 
and the $2$-bishop's move is restricted to $(0,1)$, which is $N$. 
Hence $(0,2)$ is $P$. More generally, $(0,y)$ is $N$ for all $y\ge 3$, 
since  $(0, y)\rightarrow (0,0)$ is legal if $0 < y \equiv 0,1 \pmod 4$ and 
$(0, y)\rightarrow (0,2)$ is legal if $ 2 < y \equiv 2,3 \pmod 4$.
\end{Ex}
\begin{Ex}\label{ex4}
Suppose the starting position is $(0,4)$ and the game is $2\W^3$. Then
the only bishop-type move is $(0,4)\rightarrow (0,3)$, so that the
roob-type options are $(0,0),(0,1),(0,2)$. Bob may block off 2 of
these positions, say $(0,0),(0,2)$. Then if Alice moves to
 $(0,1)$ she will loose (since she may not block off (0,0)), so suppose
rather that she moves to $(0,3)$. Than she may not block off $(0,2)$ so
Bob moves $(0,3)\rightarrow (0,2)$ and blocks off $(0,0)$. Hence
$(0,4)$ is a $P$-position.
\end{Ex}
\begin{Ex}\label{ex5}
Suppose the starting position is $(0,4)$ and the game is $2\W^{(3)}$.
Alice cannot move to $(0,0)$ or $(0,2)$. But
 $(0,1)\rightarrow (0,0)$ is a $2$-bishop-type option and
$(0,3)\rightarrow (0,0)$ is a $3$-rook-type option. This shows that
$(0,4)$ is a $P$-position.
\end{Ex}

Notice that, in comparison to
Examples \ref{ex4} and \ref{ex5}, the $P$-positions in the
Examples \ref{ex1} and \ref{ex2} are
distinct in spite the identical game constants ($m=p=2$). On the other hand,
the $P$-positions in Examples \ref{ex1} and \ref{ex3} coincide.

\begin{Ex}\label{ex6}
Suppose the game is $2\!\!\times \!\!3\W_1$, then the terminal positions 
are $(2,0)$ and $(0,4)$. On the other hand, for the game 
$2\!\!\times \!\!3\W_2$, the positions $(0,2)$ and $(4,0)$ are terminal.
Suppose now that the starting position of $2\!\!\times \!\!3\W_2$ is
$(1, 9)$. Then Alice wins by moving to $(0, 4)$. If the starting position 
is the same, but the game is $2\!\!\times \!\!3\W_1$, then Alice cannot 
move to $(0,2)$ and hence Bob wins. 
Similarly, if the starting position of $2\!\!\times \!\!3\W_0$ is
$(1, 7)$ Alice may not move to $(0, 0)$ and hence Bob wins.
\end{Ex}

\begin{figure}[h]
\centering
\includegraphics[width=1\textwidth]{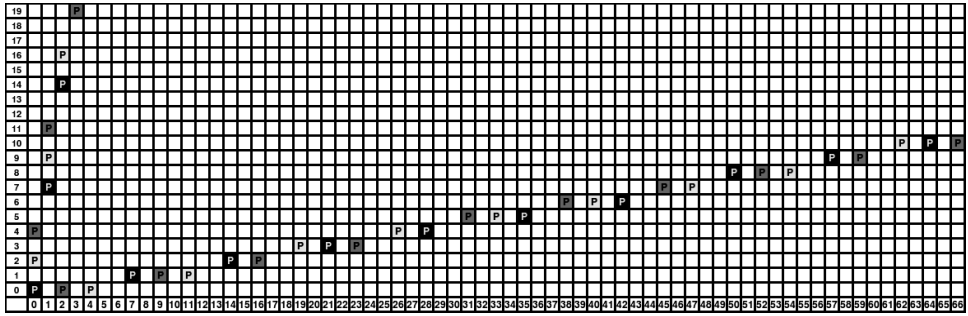}
\caption{$P$-positions of $2\W^{(3)}$, $2\W^{3}$, $2\W^{2,6}$
and $2\!\! \times \!\! 3\W$---the positions
nearest the origin such that there are precisely three positions in
each row and column and one position in every second NE-SW-diagonal.
The palest colored squares represent $P$-positions of
$2\!\!\times\!\!3\W_1$. They are of the form $(a_{3n+1}, b_{3n+1})$
or $(b_{3n+2}, a_{3n+2})$. The darkest squares, $(\{a^{2,3}_{3i},b^{2,3}_{3i}\})$,
represent the solution of $6\W$. }
\end{figure}

\begin{table}[h!]
\begin{center}
  \begin{tabular}
{| l || c | c | c | c | c | c | c | c | c | c | c | c | c | c | c | c | c|}
    \hline
    $b_n^{2,3}$\T &0&2&4&7&9&11&14&16&19&21&23&26&28&31&33&35&38 \\
    $a_n^{2,3}$\T &0&0&0&1&1&1 &2 &2 &3 &3 &3 &4 &4 &5 &5 &5 &6 \\ \hline
$b_n\!-\!a_n$\T &0&2&4&6&8&10&12&14&16&18&20&22&24&26&28&30&32\\\hline
    $n$   \T &0 &1&2&3& 4& 5& 6& 7& 8& 9& 10&11&12&13&14&15&16\\
    \hline
  \end{tabular}
\end{center}\caption{Some initial values of the Beatty pairs defined
in (\ref{A}) and (\ref{B}), here $m=2$ and $p=3$, together with the
differences of their coordinates (=$2n$).}
\end{table}

\subsection{Game theory results}
We may now state our main results. We prove them in Section 6, since
our proofs depend on some arithmetic results presented in Section 3, 4 and 5.

\begin{Thm}[Main Theorem]\label{maintheorem}
Fix $m,p\in \N$ and let $a$ and $b$ be as
in (\ref{A}) and (\ref{B}). Then
\begin{enumerate}[(i)]
\item $\mathcal{P}(m\W^p)=\{\{a_i,b_i\}\mid i\in \M\}$;
\item
\begin{enumerate}[(a)]
\item $\mathcal{P}(m\W^{(p)})=\{\{a_i,b_i\}\mid i\in \M\}$ if and
only if $\gcd(m,p)=1$;
\item $\mathcal{P}(m\W^{(m,mp)})=\{\{a_i,b_i\}\mid i\in \M\}$;
\end{enumerate}
\item
\begin{enumerate}[(a)]
\item $\mathcal{P}(m\!\!\times \!\!p\W_l)=
\{(a_{ip+l},b_{ip+l})\mid i\in \M\}\cup \{(b_{ip-l},a_{ip-l})\mid i\in \N\}$
\item $\mathcal{P}(m\!\!\times \!\!p\W)=
\{\{a_{i},b_{i}\}\mid i\in \M\}.$
\end{enumerate}
\end{enumerate}
\end{Thm}
By this result it is clear that, in terms of game complexity, the solution
of each of our games is polynomial.

\section{More on $p$-complementary Beatty sequences}
As we have seen, it is customary to represent
the solution of 'a removal game on two heaps of tokens' as a sequence of
(ordered) pairs of non-negative integers. However, often it turns
out that it is more convenient to study
the corresponding pair of sequences of non-negative integers. Sometimes, as
in Wythoff Nim, these sequences are increasing. It turns out that for our
purpose we are more interested in pairs of non-decreasing sequences.
This leads us to a certain extension of Beatty's original theorem, to
pairs of $p$-complementary Beatty sequences.

In the literature there is a proof of this theorem in
 \cite{Bry02}, where K. O'Bryant uses generating functions (a method
 adapted from \cite{BoBo93}). Here, we have chosen to
include an elementary proof, in analogy to ideas presented
in \cite{HyOs27, Fra82}. (See also the Appendix.)

\begin{Thm}[O'Bryant]\label{Beam}
Let $p\in \N$ and let $0 < \alpha < \beta $ be real numbers such
that

\begin{align}\label{abp}
\frac {1}{\alpha} + \frac {1}{\beta} = p.
\end{align}

Then $ (\lfloor i\alpha\rfloor )_{i\in \M}$ and
$ (\lfloor i\beta \rfloor )_{i\in \N}$ are $p$-complementary on $\M$ 
if and only if $\alpha , \beta $ are irrational.
\end{Thm}

\noindent {\bf Proof.} It suffices to establish that
exactly $p$ members of the set $$S =
\{0, \alpha , \beta , 2\alpha , 2\beta ,\ldots\}$$
is in the interval $[n, n+1)$ for each $n\in \M $. But for any
fixed $n$ we have
\begin{align*}
\# (S \cap [0, n])&=
\#(\{0,\alpha, 2\alpha,\ldots\}\cap [0, n]) + \#(\{\beta, 2\beta, \ldots
\}\cap[1, n])\\
&= \lfloor n/\alpha \rfloor + 1 + \lfloor n/\beta\rfloor .
\intertext{But $\alpha$ and $\beta$ are irrational if and only if,
for all $n$, }
np - 1 = n/\alpha +n/\beta -1 &< \lfloor n/\alpha \rfloor + 1 + \lfloor
n/\beta \rfloor \\ &< n/\alpha +n/\beta +1 = np + 1.
\end{align*} This gives $\lfloor n/\alpha \rfloor + 1 + \lfloor
n/\beta \rfloor = np$. Going from $n$ to $n+1$ gives the result. $\hfill\Box$\\

The following result establishes that
a pair of $p$-complementary homogeneous Beatty sequences can always 
be partitioned into $p$ complementary pairs of Beatty sequences.
For the proof we use a generalization of Beatty's
theorem in \cite{Sko57, Fra69, Bry03}.

\begin{Prop}\label{Bea2}
Let $2\le p\in \Z$. Suppose that 
$(x_i) = (\lfloor \alpha i\rfloor)_{i\in \M}$ and
$(y_i) = (\lfloor \beta i\rfloor)_{i\in \N}$ are $p$-complementary
homogeneous Beatty sequences with $0 < \alpha <\beta $.
 Then, for any fixed integer $0 \le l < p$, 
the pair of sequences $$(x_{pi+l})_{i\in \M} \text{ and } (y_{pi-l})_{i\in \N}$$ 
is complementary on $\M $.
\end{Prop}

\noindent{\bf Proof.}
Since $(x_i)$ and $(y_i)$ are non-decreasing and $p$-complementary,
we get that both $(x_{pi + l}) \text{ and } (y_{pi - l})$ are increasing.
Then, by $x_l = \min \{x_{pi + l}\}$ and $y_{p-l} = \min \{y_{pi - l}\}$, 
again, $p$-complementarity gives that
$$\max \{x_l, y_{p-l}\} > \min\{x_l, y_{p-l}\} = 0$$ 
and so, we may conclude that the integer $0$
occurs precisely once together in $(x_{pi+l}) \text{ and } (y_{pi-l})$.
Let us adapt to the terminology in \cite{Bry03}.\\

\noindent {\it Case $x_l=0$}:  Then $l < \frac{1}{\alpha}$. 
For all $n\in \N$, we denote 
$$\left \lfloor \frac{n - \gamma '}{\gamma}\right \rfloor = x_{pn + l}$$ and
$$\left \lfloor \frac{n - \eta '}{\eta}\right \rfloor = y_{pn - l}.$$
This gives $\gamma = \frac{1}{p\alpha} $, $\gamma' = -\frac{l}{p}$,
$\eta = \frac{1}{p\beta} $ and $\eta' = \frac{l}{p}$.
Then according to Fraenkel's Partition Theorem in \cite[page 5]{Bry03}, since
$\alpha$ is irrational, $(x_{pn+l}) \text{ and } (y_{pn-l})$ are complementary 
on $\N$ if and only if
\begin{enumerate}[(i)]
\item $0 < \gamma < 1$,
\item $\gamma + \eta = 1$
\item $0\le \gamma + \gamma ' \le 1$
\item $\gamma ' + \eta ' =0$ and $k\gamma + \gamma' \not\in \Z$
for $2\le k\in \Z$.
\end{enumerate}
These four items are easy to verify.\\ 
\noindent Item (i): This follows since
(\ref{abp}) together with $0 < \alpha < \beta $
is equivalent to $\frac{1}{p} < \alpha < \frac{2}{p}$.\\
\noindent Item (ii): This is immediate by (\ref{abp}).\\
\noindent Item (iii):  We have that $0\le \gamma + \gamma ' \le 1$
if and only if $0 \le \frac{1}{p\alpha} -
\frac{l}{p}\le 1$ if and only if 
$l\le  \frac{1}{\alpha} \le p + l$.\\
\noindent Item (iv): We have that
$\eta ' + \gamma '= \frac{l}{p} -  \frac{l}{p} = 0$. Since $\gamma$
is irrational and $\gamma '$ is rational the latter part holds as well.\\

\noindent {\it Case $y_{p-l}=0$}:  Then 
$p-l < \frac{1}{\beta} = p - \frac{1}{\alpha}$, so that $l > \frac{1}{\alpha}$. 
For all $n\in \N$, we denote 
$$\left \lfloor \frac{n - \gamma '}{\gamma}\right \rfloor = x_{p(n-1) + l}$$ and
$$\left \lfloor \frac{n - \eta '}{\eta}\right \rfloor = y_{p(n+1) - l}.$$
This gives $\gamma = \frac{1}{p\alpha} $, $\gamma' = \frac{p-l}{p}$,
$\eta = \frac{1}{p\beta} $ and $\eta' = \frac{l-p}{p}$. Then items 
(i), (ii) and (iv) are treated in analogy with the first case. 
For (iii), since $\alpha l >1$, we get  
$0\le \gamma + \gamma ' = \frac{1}{p\alpha} +
\frac{p-l}{p} = 1 + \frac{1-\alpha l}{\alpha p} < 1$.
\hfill $\Box$\\

We will now focus on the properties of the sequences $a$ and $b$.
The next result is central to the rest of the paper.
\begin{Lemma}\label{mainlemma}
Fix $m,p\in \N$ and let $a$ and $b$ be
as in $(\ref{A})$ and $(\ref{B})$ respectively.
Then for each $n\in \M$ we have that
\begin{enumerate}[(i)]
\item $a$ and $b$ are $p$-complementary;
\item $b_n-a_n=mn$;
\item if $p=1$, then
\begin{enumerate}
\item $a_{n+1} - a_n = 1$ and $b_{n+1} - b_n = m + 1$, or
\item $a_{n+1} - a_n = 2$ and $b_{n+1} - b_n = m + 2$;
\end{enumerate}
\item if $p > 1$, then
\begin{enumerate}
\item $a_{n + 1} - a_n = 0$ and $b_{n + 1} - b_n = m$, or
\item $a_{n + 1} - a_n = 1$ and $b_{n + 1} - b_n = m + 1$.
\end{enumerate}
\end{enumerate}

\end{Lemma}
\noindent {\bf Proof.}
Since $\phi_x$ is irrational and
$\frac{1}{\phi_x} + \frac{1}{\phi_x + x} = 1,$ case (i) is immediate from
Theorem \ref{Beam}.

For case (ii) put $\nu = \nu (m,p) = \frac{\phi_{mp}}{p} + \frac{m}{2}$
and observe that
$$b_n - a_n = \left\lfloor n\left(\nu + \frac{m}{2}\right)\right\rfloor
- \left\lfloor n(\nu - \frac{m}{2})\right\rfloor.$$ If $mn$ is even,
we are done, so suppose that $mn - 1 = 2k$, $k\in \M$. Then
$$b_n - a_n = \left\lfloor n\nu + \frac{1}{2}\right\rfloor
- \left\lfloor n\nu - \frac{1}{2}\right\rfloor + 2k=1+2k=mn.$$

For case (iii), by \cite{Fra82}, we are
done. In case $p>1$, by the triangle inequality, we get
\begin{align*}
0&<\frac{\phi_{mp}}{p}\\
&=\frac{1}{p}-\frac{m}{2}+\sqrt{\frac{m^2}{4}+\frac{1}{p^2}}\\
&<\frac{1}{p}+\frac{1}{p}\\
&\le 1, \text{ since } p\ge 2,
\end{align*}
so that we may estimate
\begin{align*}
a_{n+1} - a_n &= \left\lfloor \frac{(n+1)\phi_{mp}}{p}\right\rfloor
- \left\lfloor \frac{n\phi_{mp}}{p}\right\rfloor \in \{0,1\}.
\end{align*}
Then by (ii) we have
\begin{align*}
b_{n+1} - b_n &= a_{n+1}+m(n+1)-a_n-mn\\
&=a_{n+1}-a_n+m,
\end{align*}
so that (iv) holds. $\hfill\Box$\\

\section{A unique pair of $p$-complementary Beatty sequences}

For fixed $p$ and $m$ we now present a certain uniqueness property
for our pair of $p$-complementary Beatty sequences (in case $p=1$
see also \cite{HeLa06} for extensive generalizations).

\begin{Thm}\label{uniquethm}
Fix $m,p \in \N$. Suppose $x = (x_i)_{i\in \M}$ and
$y = (y_i)_{i\in \M}$ are non-decreasing sequences of non-negative
integers.
Then the following two items are equivalent,
\begin{enumerate}[(i)]
\item $x$ and $y_{>0}$ are $p$-complementary and,
for all $n$, $y_n - x_n = mn$;
\item for all $n$, $x_n = a_n^{m,p}$ and $y_n = b_n^{m,p}$.
\end{enumerate}
\end{Thm}
\noindent {\bf Proof.} By Lemma \ref{mainlemma} it is clear that (ii)
implies (i). Hence, it suffices to prove the other direction.

It is given that $x_0 = y_0 = a_0 = b_0 = 0$. Since $x$ is non-decreasing 
the condition $y_n - x_n = mn$ implies that $y$ is increasing.
Suppose that Lemma \ref{mainlemma} (iv) holds for a fixed $n\ge 0$,
but with $a$ exchanged for $x$ and $b$ exchanged for $y$. Then,
since $x$ and $y_{>0}$ are $p$-complementary and $y_{n+1} > x_{n+1}$,
we must have that $x_{n+1} - x_n = 0$ if
\begin{align*}
\#\{ i\mid x_i = x_n \text{ or } y_{i+1} = x_n, 0\le i\le n\} < p,
\end{align*}
and $x_{n+1} - x_n = 1$ if
\begin{align*}
\#\{ i\mid x_i = x_n \text{ or } y_{i+1} = x_n, 0\le i\le n\} = p.
\end{align*}
 But, by Lemma \ref{mainlemma}, this also holds for
the sequence $(a_i)$. In conclusion,
$y_{n+1} = x_{n+1} + m(n+1) = a_{n+1} + m(n+1) =b_{n+1}$ gives the result.
\hfill$\Box$

\section{Recurrence results}
In this section we generalize the minimal exclusive algorithm in (\ref{framex}).
Since our game rules are three-folded we will study
three different recurrences. But first we explain why 
$a$ and $b$ satisfy the '$p$-complementary equation' in (\ref{pcomp}) .

\begin{Thm}\label{recurrthm}
Fix $m,p\in \N$ and let $a$ and $b$ be as in (\ref{A}) and (\ref{B}).
For each $n\in \M$, define
$$\varphi_n=\varphi_n(m,p):= \frac{a_n+(mp-1)b_n}{m}.$$
Then, for each $n\in \N$, $\varphi_n$ is the greatest integer
such that
\begin{align}\label{recurr}
b_n-1=a_{\varphi_n}.
\end{align}
\end{Thm}

\noindent {\bf Proof.}
Notice that, for all $n$,
\begin{align}
\varphi_n&= \frac{a_n+(mp-1)b_n}{m}\notag \\
&= \frac{mpb_n-mn}{m}\notag \\
&=pb_n-n,\label{pbnn}
\end{align}
so that
\begin{align}\label{pigeon}
\varphi_{n+1}-\varphi_{n}&=pb_{n+1}-(n+1)-(pb_n-n)\notag \\
&= p(b_{n+1}-b_n)-1.
\end{align}

The proof is by induction. For the base case, notice that $b_1 = m$, $a_1 = 0$
and $\varphi_1 = (mp - 1)$. 
The only representative from
$b$ in the interval $[0, p-1]$ is $b_0 = 0$ (which we by definition 
do not count). Hence, by $a_0 = 0$ and $p$-complementarity, we get that 
$$a_{\varphi_{1}} = a_{mp-1} = m - 1 = b_1 - 1$$ and 
$$a_{\varphi_{1}+1} = a_{mp} = m = b_1.$$

Suppose that (\ref{recurr}) holds for all $i\le n$.
Then we need to show that $b_{n+1}-1=a_{\varphi_{n+1}}$
and $b_{n+1}=a_{\varphi_{n+1}+1}$.

If $a_{\varphi_{n+1}} - a_{\varphi_n} = b_{n+1}-b_n$, by $b_n - 1 = a_{\varphi_n}$
and $b_n = a_{\varphi_n+1}$, we are done, so assume that either
\begin{itemize}
\item [(A)] $a_{\varphi_{n+1}}-a_{\varphi_n} < b_{n+1}-b_n$, or
\item [(B)] $a_{\varphi_{n+1}}-a_{\varphi_n} > b_{n+1}-b_n.$
\end{itemize}

Again, by $p$-complementarity, the total number of elements 
from $a$ and $b$ in the interval
\begin{align*}
I_n &:= (a_{\varphi_n}, a_{\varphi_{n+1}}]\\
&= (a_{\varphi_n}, a_{\varphi_{n}+p(b_{n+1}-b_n)-1)}]
\end{align*}
is $R_n:=p(a_{\varphi_{n+1}}-a_{\varphi_n})$, and where the equality
is by (\ref{pigeon}). By assumption,
$a_{\varphi_n+1}\in I_n$ so that we have at least
$p(b_{n+1}-b_n)-1$ representatives 
from $a$ in $I_n$. But also $b_n=a_{\varphi_n}+1\in I_n$ so that
altogether we have at least $p(b_{n+1}-b_n)$ representatives in $I_n$.
Hence
\begin{align*}
p(b_{n+1}-b_n)&\le R_n \\
&= p(a_{\varphi_{n+1}}-a_{\varphi_n})
\end{align*}
which rules out case (A).

Notice that case (B) implies that $b_{n+1}$ lies in $I_n$ so that
$a_{\varphi_{n}+1} = b_n < b_{n+1} \le a_{\varphi_{n+1}}$. Since both $b_n$
and $b_{n+1}$ lie in $I_n$, the total number of representatives in $I_n$ 
is 
\begin{align}
p(a_{\varphi_{n+1}} - a_{\varphi_{n}})
&\le 2 + \varphi_{n+1}-\varphi_{n}\notag\\
&=p(b_{n+1} - b_n) + 1.\label{tag}
\end{align}
In case $p>1$, since $a$ and $b$ are integer sequences, 
we are done, so suppose $p=1$. Then, in fact, by complementarity, we must have 
$a_{\varphi_{n+1}}<b_n<b_{n+1}<a_{\varphi_{n+1}},$ contradicting (\ref{tag}). 
$\hfill\Box$


\begin{Rem}
For arbitrary $m>0$ and $p=1$ it is well-known that $a$ and $b$ solve
$x_{y_n}=x_n+y_n$. This complementary equation is studied in for example
\cite{Conn59, FrKi94, Kim07}. However, we have not been able to find any
references for the complementary equation $y_n - 1 = x_{y_n-n}$.  
By (\ref{pbnn}), for the cases $p=1$, this equation is also resolved 
by $a$ and $b$.
\end{Rem}

A \emph{multiset} (or a sequence) $X$ may
be represented as (another) sequence of non-negative
integers $\xi = (\xi^i)_{i\in \M}$,
where, for each $i\in \M$, $\xi^i=\xi^i(X)$ counts the
number of occurrences of $i$ in $X$.
For a positive integer $p$, let $\mex^p \xi $
denote the least non-negative integer $i\in X$ such that $\xi^i < p.$

\begin{Prop}\label{mexprop}
Let $m, p\in \N$.  Then the definitions of the sequences $x$ and $y$ in 
(i), (ii) and (iii) are equivalent. In fact,
for each $n\in \M$, we have that $x_n = a_n^{m,p}$ and $y_n = b_n^{m,p}$. 
\begin{enumerate}[(i)]
\item For $n\ge 0$,
\begin{align*}
x_n &=\mex^p\xi_n,\\
\intertext{where $\xi_{n}$ is the multiset, where for each $i\in \M,$}
\xi^i_{n} &= \#\{j \mid i=x_j  \text{ or } i=y_j, 0\le j< n\},\\
y_n&=x_n+mn.
\end{align*}
\item For $n\ge 0$,
\begin{align*}
x_n &= \mex \{\nu_i^n,\mu_i^n \mid  0\le i<n\},\text{ where }\\
\nu_i^n&=x_i \text{ if } n \equiv i\!\!\! \pmod {p},
\text{ else } \nu_i^n=\infty, \\
\mu_i^n&=y_i \text{ if } n \equiv -i\!\!\! \pmod {p},
\text{ else }\mu_i^n=\infty ,\\
y_n &= x_n + mn.
\end{align*}
\item For $n\ge 0$,
\begin{align*}
x_{pn} &= \mex \{ x_{pi}, y_{pi} \mid 0\le i < n \},\\
y_{pn} &= x_{pn} + mpn,\intertext{ and for each integer $0 < l < p$,}
x_{pn+l} &= \mex \{ x_{pi+l}, y_{p(i+1)-l} \mid 0\le i < n \},\\
y_{pn+l} &= x_{pn+l} + m(pn + l).
\end{align*}
\end{enumerate}
 
\end{Prop}

{\bf Proof.}
For $p=1$ it is a straightforward task to check that 
each recurrence is equivalent to (\ref{framex}).
Hence, let $p > 1$. Observe that in (i), by definition, $x$ and $y$ 
are non-decreasing, $p$-complementary and, for all $n$, 
\begin{align}\label{yxmn}
y_n = x_n + mn.
\end{align}
Hence, for this case, Theorem \ref{uniquethm} gives the result.

Let us now study the definitions of $x$ and $y$ in (ii). 
For $z\in\Z$, let $\overline z$ denote the congruence 
class of $z$ modulo $p$. Here, it is not immediately clear that the sequences 
are non-decreasing. Neither is it obvious that they are $p$-complementary. 
But, at least we have that, for each $n\in \M$, (\ref{yxmn}) holds.

Hence, if (ii) fails (by $a_0 = x_0$) there has to exist a least 
index $n'\in \N$ such that $a_{n'}\ne x_{n'}$. But notice that 
$0\le n < p$ implies $\nu_i^n = \mu_i^n = \infty ,$ for all $0\le i < n$, which 
in its turn implies $a_n = x_n = 0$. This gives $n'\ge p$. We have two 
cases to consider: 
\begin{enumerate}[(a)]
\item $r := x_{n'} < a_{n'}:$  By Theorem \ref{recurrthm}
there are two cases to consider.
\begin{itemize}
\item [Case 1:] There is an $i\ge 0$ such that $\varphi (i)+p-1<n'$
and $$y_i=x_{\varphi(i)+1}=x_{\varphi(i)+2}=\ldots =x_{\varphi(i)+p-1}=r.$$
But then, by
\begin{align}\label{harder}
\{\;\overline{-i},\; \overline{-i+1},\; \ldots
, \;\overline{-i+p-1}\; \}=
\{\; \overline 0,\; \overline 1,\; \ldots ,\; \overline{p-1}\; \}
\end{align}
and 
\begin{align}\label{modm}
\varphi_n = pb_n - n \equiv -n \pmod p,
\end{align}
there is a $j\in \{i,\varphi(i)+1,\ldots , \varphi(i)+p-1\}$ such
that either  $n'\equiv j\pmod p$ and
$j\in\{\varphi(i)+1,\ldots \varphi(i)+p-1\}$  which
implies $\nu_j^{n'}=r$, or  $n'\equiv -j\pmod p$ and $j=i$ which implies
$\mu_j^{n'}=r$. In either case the
choice of $x_{n'}=r$ contradicts the definition of $\mex$.
\item [Case 2:] There is an $i\ge 0$ such that $i+p-1<n'$
and $$r=x_i=x_{i+1}=x_{i+2}=\ldots =x_{i+p-1}.$$
This case is similar but simpler, since for this case we rather use that
\begin{align}\label{simpler}
\{\;\overline{i},\; \overline{i+1},\; \ldots
, \;\overline{i+p-1}\; \}=
\{\; \overline 0,\; \overline 1,\; \ldots ,\; \overline{p-1}\; \}
\end{align}
\end{itemize}
\item $r:=a_{n'}<x_{n'}:$ Then our $\mex$-algorithm has refused $r$ as the
choice for $x_{n'}$. But
then there must be an index $0 \le j < n'$ such that either $\nu_j^{n'} = r$
or $\mu_j^{n'} = r$. Hence, we get to consider two cases.
\begin{itemize}
\item[Case 1:] $\overline{j}\; =\overline{n'}$ and $r = x_j$. On
the one hand, there is a $k\in \N$ such that $kp + j = n'$ On the other
hand,  there is a greatest $k'\in \N$ such that
$a_{n'-k'} = a_{n'-k'+1} = \ldots = a_{n'}$ and by $p$-complementarity
$0\le k' < p$. But then, since $n'-k' > n'- kp = j$,
we get $a_j < r = x_j$, which contradicts the minimality of $n'$.
\item[Case 2:] $\overline{-j} \; =\overline{n'}$ and $r = y_j$. Then, 
by Theorem \ref{recurrthm}, $\varphi_j + 1$ is the least index such that
$a_{\varphi_j+1} = a_{n'}$. Then, by minimality of $n'$, $a_j = x_j$ gives 
$b_j = y_j$ so that $a_{n'} = b_j$. For this case, $p$-complementarity gives
$n' - (\varphi_j + 1) + 1 \le p - 1$.  Then $0 < k' := n'- \varphi_j < p$ 
and so $$\overline{-j + k'} = \; \overline{\varphi(j) + k'}
=\; \overline{n'} = \; \overline{-j},$$ which is nonsense.
\end{itemize}
\end{enumerate}

For case (iii), suppose that there is a least index $n' \ge p$ such
that $a_{n'} \ne x_{n'}.$ (The case $n' < p$ may be ruled out 
in analogy with (ii).) Then, there exist unique integers, $0 < t$ and
$0 \le l < p$, such that $tp + l = n'$.

Suppose that $r := a_{n'} < x_{n'}$. Then, since
the $\mex$-algorithm did not choose
$x_{n'} = r$, there must be an index $0 \le t' < t$ such
that either $x_{t'p+l} = r$ or $y_{(t'+1)p-l} = r$. But then, by assumption, 
either $a_{t'p+l} = x_{t'p+l} = a_{n'}$ or 
$b_{(t'+1)p-l} = y_{(t'+1)p-l} = a_{n'}$. But, by Proposition \ref{Bea2} $a$ 
and $b$ are complementary, so either case is ridiculous. 

Hence, assume $r := a_{n'} > x_{n'}$.
Then again, by Proposition \ref{Bea2}, there is an index $0 \le t' < t$
such that either $a_{t'p+l} = x_{n'}$ or $b_{(t'+1)p-l} = x_{n'}$. But, again, 
by minimality of $n'$, this contradicts the $\mex$-algorithm's choice 
of $x_{n'} < a_{n'}$. $\hfill \Box$

\section{The games final section}
By comparing Definitions \ref{moves}, \ref{roob} and \ref{games} 
with the results in Proposition \ref{mexprop} one could, 
already at this point, claim that 
Theorem \ref{maintheorem} holds. Namely, the allowed moves in the 
the respective games defined the recurrence in the 
three Minimal EXclusive algorithms, which in turn precisely 
determined the $P$-positions of the respective games. But, of course, 
the paper would remain incomplete without an explicit 
\emph{game-theoretical} proof.\\

\noindent{\bf Proof of Theorem \ref{maintheorem}.}
For $p = 1$, the games have identical
rules. This case has been established in \cite{Fra82}.
The case $m = 1$ has been studied in \cite{Con59} for games of form (ii).
(and implicitly for $1\!\!\times \!\!p\W_l$).

For the rest of the proof assume that $p > 1$. 
For each game we need to prove that, if $(x, y)$
\begin{enumerate}[(A)]
\item is of the form $\{a_i, b_i\}$, then none of its options is;
\item is not of the form $\{a_i, b_i\}$, then it has an option
of this form.
\end{enumerate}
(We will need a slightly different notation for Case (iii) below.)
By symmetry, we may assume that $0\le x\le y$. Clearly, any final 
position satisfies (A) but not (B).

\begin{itemize}
\item [Game (i):] We need to prove that 
$\mathcal{P}(m\W^p)=\{\{a_i, b_i\}\mid i\in \M\}$. Suppose 
$(x, y)=(a_i, b_i)$ for some $i\in \M$. By Lemma \ref{mainlemma} (i) and (ii),
$a$ and $b_{>}$ are $p$-complementary and $b_i - b_j\ge m$ for all $j < i$.
Then any roob-type option of the form $\{a_j, b_j\}$ may be blocked off, unless
perhaps $a_{j} < a_i$ and $b_{j} = b_i$ for some $j < i$.
But this is ridiculous since $b$ is strictly increasing. By
Lemma \ref{mainlemma} (ii) we get that,
for $j < i$, $b_i - a_i\pm (b_j - a_j)\ge m$.
Then an $m$-bishop cannot move
$(x, y)\rightarrow \{a_j, b_j\}$, This proves (A).

For (B), since $p\ge 2$,
we may assume $x = a_i$, for some $i$, but $y\ne b_i$.
Then, by Lemma \ref{mainlemma} (iv): (*) There
exists a $j<i$ such that an $m$-bishop can move $(x,y)\rightarrow (a_j,b_j)$
unless $y-x-(b_j-a_j)\ge m$ for all $j$ such that $a_j\le x$.
Then, for all $j$ such that $a_j = x$, we have that 
$y\ge x+m(j+1)> b_j$. But then, by Lemma \ref{mainlemma} (i), there are $p$ 
options of $(x,y)$ of the form $\{a_i, bi\}$. By the rule of game, 
they cannot all be blocked off.\\

\item [Game (iia):]  We are going to prove that 
$\mathcal{P}(m\W^{(p)}) = \{\{a_i,b_i\}\mid i\in \M\}$ if and
only if $\gcd(m,p) = 1$. Let us first explain the 'only if' direction. 
Denote with $\gamma = \gcd(m, p)$, $p' = \frac{p}{\gamma }$ and 
$m' = \frac{m}{\gamma} $. Then, clearly,
the positions of the form $(0, mi)$, where $0\le i < p' $,
are $P$-positions of $m\W^{(p)}$. Now,
$(0,mp')$ is an $N$-position because $m'p=mp'$ implies that 
$(0,mp')\rightarrow (0,0)$ is an option. 
But, by definition, $b_{p'} = mp'$ if and only if $p' < p$ if and only if
$\gamma >1$. Hence $\gcd(m, p) = 1$ is a necessary requirement.

For this game, the options of the $m$-bishop
are identical to those in (i). Hence, let us analyze the $p$-rook.

For (A), suppose that $(x, y) = (a_i, b_i)$ for some $i\in \M$ 
but that, for a contradiction, that a $p$-rook can move to $\{a_j, b_j\}$. 
Then, since $b$ is strictly increasing, there is a $0\le j < i$, such 
that either $b_i \equiv b_j \pmod p$
and $a_i = a_j$, or $b_i \equiv a_j \pmod p$ and $a_i = b_j$. But then,
for the first case (using the same notation as in Section 5), 
since $$\overline {mj} = \overline{b_j - a_j}
= \overline{b_i - a_i} = \overline {mi}$$ and $\gcd(m,p) = 1$ we
must have $\overline {j} = \overline {i}$. This is ridiculous, since
by $p$-complementarity and $a$ non-decreasing we have $0 < i - j < p$.
For the second case, by Theorem \ref{recurrthm}, we have that
$$\overline{-mj} = \overline{a_j - b_j}
= \overline{b_i - a_i} = \overline {mi} = \overline {m(\varphi (j) + t)}
= \overline{m(-j + t)},$$ 
for some $t\in \{1, \ldots , p - 1\}$. This implies 
$\overline 0 = \overline {mt}$ but then again $\gcd(m, p) = 1$ gives
a contradiction.

For (B), we follow the ideas in the second part of Case (i)
up until (*). Then, for this game, we rather need to show
that there is a $j < i$ such that $y\equiv b_j \pmod p$ and $a_j = x$
or $y\equiv a_j \pmod p$ and $b_j = x$. But this follows directly from
the proof of Proposition \ref{mexprop} (ii)(a).\\

\item  [Game (iib):] We are now going to show that 
$\mathcal{P}(m\W^{(m,mp)}) = \{\{a_i,b_i\}\mid i\in \M\}$. 
For (A), suppose $(x, y) = (a_i, b_i)$ for some $i \in \M$ but
that there is a $j < i$ such that 
the $(m, mp)$-rook can move to $\{a_j, b_j\}$.
Then, we have two cases:
\begin{itemize}
\item [Case 1:] $b_i\equiv b_j - r \pmod {mp}$ and $a_i = a_j$, for some
$r\in \{0, 1, \ldots , m - 1\}$. Then
$b_i - a_i\equiv b_j - a_j - r \pmod {mp}$ so that
$mi \equiv mj - r \pmod {mp}$ and so
$m(i - j)\equiv - r \pmod {mp}$. But this forces $r = 0$ and 
$i - j\equiv 0 \pmod p$ which is impossible since Lemma \ref{mainlemma} (i) 
and (iv) imply $i - j\in \{1, 2, \ldots , p - 1\}$.
\item [Case 2:] $b_i \equiv a_j - r \pmod {mp}$ and $a_i = b_j$, for some
$r\in \{0,1,\ldots , m-1\}$. Then $b_i-a_i\equiv a_j-b_j-r \pmod {mp}$ so that
$mi\equiv -mj-r \pmod {mp}$ and so $m(i + j)\equiv - r \pmod {mp}$. 
By Theorem \ref{recurrthm} we have that $i = \varphi(j) + s$ for some 
$s\in \{1, 2, \ldots , p - 1\}$. Further, by
(\ref{modm}), we have $\varphi(j)\equiv -j \pmod p$, so that
$m(\varphi(j)+s+j) = ms \equiv - r \pmod {mp}$. Once again we
have reached a contradiction.
\end{itemize}
For (B), in analogy with (*), it suffices to study
the $(m, mp)$-rook's options where $y$ is such that
$y - x - (b_j - a_j)\ge m$ for all $j$ such that $a_j\le x = a_i$.
Hence, we need to show that there are a $j$ and
an $r\in \{0, 1, \ldots m - 1\}$ such that
$$y\equiv b_j - r \!\!\pmod {mp}\;\;\text{ and }\;\; a_j = x,$$ or
$$y\equiv a_j - r \!\!\pmod {mp}\;\; \text{ and }\;\; b_j = x.$$ Clearly,
we may choose $r$ such that $y - x + r \equiv 0\pmod m$.
Then, for all $j$, we get $ms := y - x + r\equiv \pm(b_j - a_j)\pmod m$. 
Hence, it suffices to find a specific $j$ such that
$$j = \frac{b_j - a_j}{m} \equiv s\!\!\pmod p \;\; \text{ and }\;\; a_j = x,$$
or 
$$-j = \frac{a_j - b_j}{m}\equiv s\!\!\pmod p \;\; \text{ and }\;\; b_j = x.$$
But then, by (\ref{harder}) or (\ref{simpler}), we are done.\\

\item [Game (iiia):] We are now going to show that 
$\mathcal{P}(m\!\!\times \!\!p\W_l) = 
\{(a_{ip+l},b_{ip+l})\mid i\in \M\}\cup \{(b_{ip-l},a_{ip-l})\mid i\in \N\}$. 
We may assume that $l > 0$. We have already seen that
 $(a'_i) := (a_{pi+l})_{i\ge 0}$ and $(b'_i) := (b_{p(i+1)-l})_{i\ge 0}$
are complementary. Our proof will be a straightforward extension of
those in \cite{Fra82} (which deals with the case $l = 0$)
and \cite{Con59} (which implicitly deals with the case $m = 0$).
Observe that $a'_0 = a_l = 0$ and $b'_0 = b_{p-l} = m(p - l)$.

For (A), let $(x, y) = (a'_i, b'_i)$.  In case $i = 0$ (by Definition
\ref{games} (iiia)), the Queen has no options at all, so assume $i > 0$.
Proposition \ref{mexprop} (iii) gives
that $b'_i - a'_i\pm (b'_j - a'_j) \ge mp$ for all $0 \le j < i$.
Then the $mp$-bishop cannot move $(x, y)\rightarrow (a'_j, b'_j)$ for
any $0 \le j < i$. Since $a'$ and $b'$ are complementary there is no
rook-type option $(a'_i, b'_i)\rightarrow \{a'_j, b'_j\}$.

For (B), we adjust the statement (*) accordingly: 
Suppose $x = a'_i$. By Proposition \ref{mexprop} (iii): 
If the $mp$-bishop cannot move to $(a'_j, b'_j)$ for any $j < i$ we get that
either $i = 0$ or $y - x - (b'_j - a'_j) \ge mp$ for all $j < i$.
If $i = 0$ there is a rook-type option to $(a'_0, b'_0)$
(we may assume here that $y > b'_0$), so suppose $i > 0$.
But then, by Proposition \ref{mexprop} (iii), 
we get $y\ge b'_j + mp + x - a'_j\ge b'_i + a'_i - a'_j > b'_i$. Hence, 
for this case, the rook-type move $(x, y)\rightarrow (a'_i, b'_i)$ suffices.
Suppose on the other hand that $x = b'_i$ with $i \ge 0$.
Then, since $y \ge x = b'_i > a'_i$, the desired rook-type move is 
$(x, y)\rightarrow (b'_i, a'_i)$.\\

\item [Game (iiib):] It only remains to demonstrate that 
$\mathcal{P}(m\!\!\times \!\!p\W) = \{\{a_{i},b_{i}\}\mid i\in \M\}.$ 
Suppose that the starting position is
$(a_i, b_i)$. Then $i = pj + l'$ for some (unique) pair $j\in \M$ 
and $0 \le l' < p$. The second player should choose $l = l'$. If, on the 
other hand, the starting position is
$(b_i, a_i)$. Then $i = pj - l'$ for some (unique) pair $j \in \N$ and 
$0 < l' \le p$. The second player should choose $l = p - l'$. In either case, 
by (iiia), there is no option of the form $(a'_i, b'_i)$.

If the starting position $(x, y)$ is not of the form $\{a_i, b_i\}$, again, 
by (iiia), for any choice of $0 \le l < p$, there is a 
move $(x, y)\rightarrow \{a'_i, b'_i\}$ for some $i\ge 0$.
\end{itemize}
\hfill $\Box$

\section{Questions}
Can one find a polynomial time solution of $m\W^{(l,p)}$ for some
integers $l \ge 0$, $m > 0$ and $p > 0$ whenever
\begin{itemize}
\item $\gcd(m, p) \ne 1$ and $l = 0$, or
\item $0 < l \neq m$ or $m \nmid p$?
\end{itemize}
If this turns out to be complicated, can one at least say something about its
asymptotic behavior?

Denote the solution
of $m\W^{(l,p)}$  with $\{\{ c_i^{(l,m,p)},d_i^{(l,m,p)}\}\}_{i\in \M}$.
Let us finish off with two tables of the initial $P$-positions of such games.
\begin{table}[ht!]
\begin{center}
  \begin{tabular}
{| l || c | c | c | c | c | c | c | c | c | c | c | c | c | c | c | c | c|}
    \hline
    $d_n^{(0,2,2)}$\T &0&3&6&9& 12&15&19&22&25&28&31&34&37&40&43&46&49 \\
    $c_n^{(0,2,2)}$\T &0&0&1&1& 2& 2& 3& 4& 4& 5& 5& 6& 7& 7& 8& 8& 9\\ \hline
    $d_n-c_n$\T   &0&3&5&8&10& 13&16&18&21&23&26&28&30&33&35&38&40\\\hline
    $n$   \T &0 &1&2&3& 4& 5& 6& 7& 8& 9& 10&11&12&13&14&15&16\\
    \hline
  \end{tabular}
\end{center}\caption{The first few $P$-positions of $2\W^2$ together
with the respective differences of their coordinates.}
\end{table}
\begin{table}[ht!]
\begin{center}
  \begin{tabular}
{| l || c | c | c | c | c | c | c | c | c | c | c | c | c | c | c | c| c|}
    \hline
    $d_n^{(1,2,3)}$\T &0 &2 &5 &7 &11&14&16&19&21&26&29&31&36&39&41&44&46 \\
    $c_n^{(1,2,3)}$\T &0 &0 &1 &1 &2 &3 &3 &4 &4 &5 &6 &6 &7 &8 &8 &9 &9\\\hline
    $d_n-c_n$\T     &0 &2 &4 &6 &9 &11&13&15&17&21&23&25&29&31&33&35&37\\\hline
    $n$   \T      &0&1&2&3& 4& 5& 6& 7& 8& 9& 10&11&12&13&14&15&16\\
    \hline
  \end{tabular}
\end{center}\caption{The first few $P$-positions of $2\W^{(1,3)}$. Notice
that (as in Table 3) the successive differences of their coordinates are not in
arithmetic progression.}
\end{table}

From these tables one may conclude that: The infinite arithmetic
progressions of the sequences
$$(b^{m,p}_i~-~a^{m,p}_i)_{i\in \M}=(mi)_{i\in \M}$$
(see also Table 2) are not in general seen among the sequences
$$( d_i^{(l,m,p)}-c_i^{(l,m,p)})_{i\in \M}.$$
We believe  that the latter sequence is an arithmetic progression
if and only if none of the items in our above question
is satisfied. We also believe that, for arbitrary constants,
$( c_i^{(l,m,p)})_{i\in \M}$ and $(d_i^{(l,m,p)})_{i\in \N}$ are $p$-complementary. 
But the solution of these questions are left for some future work.
\begin{Rem}
We may also define generalizations of $m\W^p$ and $m\times p\W_l$:

Fix $l\in \N$. Let $m\W^p_l$ be as $m\W^p$ but
where the player may only block off \emph{$l$-roob-type options}
(recall, non-$l$-bishop options). Otherwise,
the Queen moves as the $m$-bishop or the rook. Then $m\W^p_m = m\W^p$.
On the other hand $m\W^p_1$ is the blocking variation of $m$-Wythoff Nim 
where the previous player may block off \emph{any} $p-1$ rook-type options.

Let $u,v\in \N$ and let $m\times p\W_{u,v}$ be as  $m\times p\W_l$, but
the removed (lower left) rectangle has base $u$ and height $v$. Then
for this game
the final positions are $(u, 0)$ and $(0, v)$. If $l > 0$, $u = ml$ 
and $v = m(p-l)$
we get $m\times p\W_{lm,m(p-l)} = m\times p\W_l$. Some of these games 
are identical to mis\`ere versions of Wythoff Nim, see \cite{Fra84}.  

One may ask questions in analogy to the above for these variations.
For example, we have found a minimal
exclusive algorithm satisfying $\mathcal{P}(m\W^p_1)$ which is
related to a polynomial time construction in \cite{Fra98}. Is there an analog
polynomial time construction for $\mathcal{P}(m\W^p_1)$? Another question
is if any of these further generalized games coincide via identical set
of $P$-positions?
\end{Rem}


\noindent {\bf Acknowledgments.} I would like to thank
Aviezri Fraenkel for providing two references that motivated
generalizations of the games (in the previous version of this paper)
to their current form and of course for the nice Appendix. 
I would also like to thank Peter Hegarty for giving valuable feedback 
during the earlier part of this work,
Niklas Eriksen for composing parts of the caption for the figures, 
Johan W\"astlund for inspiring discussions about games with a 
blocking maneuver and Kevin O'Bryant for some valuable email correspondence. 
At last I would like to thank the anonymous referee
for several suggestions that helped to improve this paper.

\appendix
The following discussion, provided by Aviezri Fraenkel, 
gives a detailed analysis of '$p$-complementarity'/'$p$-fold
complementarity' and homogeneous Beatty sequences:\\

\noindent{\bf Definition 1.}
{\rm Let $p\in\N $. The multisets $S$, $T$ of positive integers are 1-{\it upper $p$-fold complementary\/}, for short: $p$-{\it fold complementary}, if $S\cup T= p\times\N $.\\}

If the multisets $S$, $T$ satisfy Definition 1 and have irrational densities $\alpha^{-1}$, $\beta^{-1}$,  say $\alpha\le\beta$, then a necessary condition for $p$-fold complementarity is  $\alpha^{-1}+\beta^{-1}=p$. Thus $a:=\beta-\alpha>0$. Then
$\alpha=(2-ap+\sqrt{a^2p^2+4})/2p,$
so $p^{-1}<\alpha<2p^{-1}$. Then $1/\beta=p-1/\alpha$, so $\beta>2/p$. 

Let $M=\lfloor 1/\alpha\rfloor+1$, $N=\lfloor 1/\beta\rfloor+1$. Notice that $\alpha (M-1) < 1 < \alpha M$, $\beta (N-1) < 1 < \beta N$. From now on we let $S=\{{\lfloor n\alpha\rfloor}\}_{n\ge M}$, $T=\{{\lfloor n\beta\rfloor}\}_{n\ge N}$.\\

\noindent{\bf Theorem 1. }{\rm The multisets $S$, $T$ 
are $p$-fold complementary.}\\

\noindent{\bf Proof.}
For any $k\in\N$, since $\alpha$ is irrational, the number of terms less than $k$ in $S\cup T$ is 
\begin{align*}
\lfloor k/\alpha\rfloor-(M-1)+\lfloor k/\beta\rfloor-(N-1)&=
\lfloor k/\alpha\rfloor+\lfloor k(p-\alpha^{-1})\rfloor-M-N+2\\
&= kp+\lfloor k/\alpha\rfloor+\lfloor -k/\alpha\rfloor-M-N+2\\
&= kp-M-N+1.
\end{align*}
 Similarly, $S\cup T$ contains $(k+1)p-M-N+1$ terms $<k+1$. Hence there are exactly $p$ terms $<k+1$ but not $<k$. They are the terms $k$ with multiplicity $p$.\hfill $\Box$\\

{\bf Remarks.}
{\rm (i)~$\lfloor 1/\alpha\rfloor=p+\lfloor -1/\beta\rfloor=p-1-\lfloor 1/\beta\rfloor=p-N$. Hence $M=\lfloor 1/\alpha\rfloor+1=p-N+1$.

(ii)~Clearly $(\lfloor (M-1)\alpha\rfloor, \lfloor (N-1)\beta\rfloor)=(0,0)$. Since $\alpha<\beta$, we have $N\le M$. 
Hence, for all $N\le n < M$, we have that $(\lfloor n\alpha\rfloor, \lfloor n\beta\rfloor) = (0,\lfloor n\beta\rfloor)$ where $\lfloor n\beta \rfloor > 0$. Thus there are precisely $M - N$ couples $(\lfloor n\alpha\rfloor, \lfloor n\beta\rfloor)$ with $\lfloor n\alpha\rfloor = 0$ and $\lfloor n\beta\rfloor > 0$, containing $M - N$ 0s.
Thus, for $0\le n < M$, there are precisely $M$ couples $(\lfloor n\alpha\rfloor, \lfloor n\beta\rfloor)$ with $\lfloor n\alpha\rfloor = 0$ and $\lfloor n\beta\rfloor \ge 0$, containing, in total, $M + N = p + 1$ 0s.

(iii) The proof is a straightforward generalization to $p\ge 1$ of a proof included in an editorial comment to \cite{AMM} stating: ``...The result is so astonishing and yet easily proved that we include a short proof for the reader's pleasure." This is then followed by the above proof for the special case $p=1$, which is itself a slight simplification of the proof given in \cite{Fra82}.}\\


\begin{thebibliography}{NOORS}
\bibitem[AlNoWo07]{AlNoWo07} M. H. Albert, R. J. Nowakowski, D. Wolfe
\emph{Lessons in Play: In Introduction to Combinatorial Game Theory}.
A K Peters Ltd.(2007).
\bibitem[AMM]{AMM}  Solution~II, Problem 11365, \emph{Amer. Math. Monthly}, 
(April 2010), p.~376,
\bibitem[Bea26]{Bea26} S. Beatty, Problem 3173, \emph{Amer. Math. Monthly},
{\bf 33} (1926) 159.
\bibitem[BeCoGu82]{BeCoGu82} E. R. Berlekamp, J. H. Conway, R.K. Guy,
\emph{Winning ways}, {\bf 1-2} Academic
Press, London (1982). Second edition, {\bf 1-4}.
A. K. Peters, Wellesley/MA (2001/03/03/04).
\bibitem[BoBo93]{BoBo93} J. M. Borwein and P. B. Borwein, On the generating
function of the integer part: $[\alpha n + \gamma ]$,
\emph{J. Number Theory} 43 (1993), pp. 293-318.
\bibitem[BoFr73]{BoFr73} I. Borosh, A.S. Fraenkel,
A Generalization of Wythoff's Game,
\emph{Jour. of Comb. Theory} (A) {\bf 15} (1973) 175-191.
\bibitem[BoFr84]{BoFr84} M. Boshernitzan and A. S. Fraenkel,
A linear algorithm for nonhomogeneous spectra of numbers,
\emph{J. Algorithms}, {\bf 5}, no. 2, pp. 187-198, 1984.
\bibitem[Bou02]{Bou02} C.L. Bouton, Nim, a game with a complete mathematical
theory, \emph{The Annals of Math. Princeton} (2) {\bf 3} (1902), 35-39.
\bibitem[Bry02]{Bry02} K. O'Bryant, A Generating Function Technique
for Beatty Sequences and Other Step Sequences,
\emph{J. Number Theory} 94, 299--319 (2002).
\bibitem[Bry03]{Bry03} K. O'Bryant, Fraenkel's Partition and Brown's
Decomposition \emph{Integers}, {\bf 3} (2003), A11, 17 pp.
\bibitem[Con59]{Con59} I.G. Connell, A generalization of
Wythoff's game \emph{Can. Math. Bull.} {\bf 2} no. 3 (1959), 181-190.
\bibitem[Conn59]{Conn59} I.G. Connell, Some properties of Beatty sequences I
\emph{Can. Math. Bull.} {\bf 2} no. 3 (1959), 190-197.
\bibitem[Con76]{Con76} J. H. Conway: \emph{On numbers and games},
Academic Press, London (1976). Second edition, A. K. Peters,
Wellesley/MA (2001).
\bibitem[DuGr08]{DuGr08} E.Duch\^ene, S. Gravier, Geometrical Extensions
of Wythoff's Game, to appear in \emph{Discrete Math} (2008).
\bibitem[Fra69]{Fra69} A.S. Fraenkel, The bracket function and complementary
sets of integers, \emph{Canad. J. Math.} {\bf 21} (1969), 6-27.
\bibitem[Fra73]{Fra73} A.S. Fraenkel, Complementing and exactly
covering sequences, J. Comb. Theory (Ser A), {\bf 14} (1973) 8-20.
\bibitem[Fra82]{Fra82} A.S. Fraenkel, How to beat your Wythoff games'
opponent on three fronts, \emph{Amer. Math. Monthly} {\bf 89} (1982) 353-361.
\bibitem[Fra84]{Fra84} A.S. Fraenkel, \emph{Wythoff games, continued fractions, 
cedar trees and Fibonacci searches}, Theoret. Comput. Sci. 29 (1984) 49-73.
\bibitem[Fra98]{Fra98} A.S. Fraenkel, Heap Games, Numeration Systems and
Sequences. \emph{Ann. of Comb.}, {\bf 2} (1998) 197-210.
\bibitem[Fra04]{Fra04} A.S. Fraenkel, Complexity, appeal and challenges
of combinatorial games. \emph{Theoret. Comp. Sci.}, {\bf 313} (2004) 393-415.
\bibitem[FrPe09]{FrPe09} A.S. Fraenkel, Udi Peled, Harnessing the
Unwieldy MEX Function, preprint, http://www.wisdom.weizmann.ac.il/~fraenkel/Papers/ \\Harnessing.The.Unwieldy.MEX.Function\_2.pdf.
\bibitem[GaSt04]{GaSt04} H. Gavel and P. Strimling, Nim with a Modular
Muller Twist, \emph{Integers: Electr. Jour. Comb. Numb. Theo.} {\bf 4} (2004).
\bibitem[Had]{Had} U. Hadad, Msc Thesis, Polynomializing Seemingly
Hard Sequences Using Surrogate Sequences,
\emph{Fac. of Math. Weiz. In. of Sci.}, (2008).
\bibitem[HeLa06]{HeLa06} P. Hegarty and U. Larsson, Permutations of the
natural numbers with prescribed difference 
multisets, \emph{Integers} {\bf 6} (2006), Paper A3, 25pp.
\bibitem[HoRe]{HoRe} A. Holshouser and H. Reiter, Three Pile Nim
with Move Blocking, http://citeseer.ist.psu.edu/470020.html.
\bibitem[HyOs27]{HyOs27} A. Ostrowski and J. Hyslop, Solution to Problem 3177,
\emph{Amer. Math. Monthly}, {\bf 34} (1927), 159-160.
\bibitem[FrKi94]{FrKi94} A.S. Fraenkel and C. Kimberling, Generalised
Wythoff arrays, shuffles
and interspersions, {\em Discrete Math.} {\bf 126} (1994), 137-149.
\bibitem[Kim95]{Kim95} C. Kimberling, Stolarsky interspersions,
{\em Ars Combinatoria} {\bf 39} (1995), 129-138.
\bibitem[Kim07]{Kim07} C. Kimberling, Complementary equations,
\emph{J. Integer Sequences} {\bf 10} (2007), Article 07.1.4.
\bibitem[Kim08]{Kim08} C. Kimberling, Complementary equations and
  Wythoff sequences, \emph{J. Integer Sequences} {\bf 11} (2008),
Article 08.3.3.
\bibitem[Lar09]{Lar09} U. Larsson, 2-pile Nim with a Restricted
Number of Move-size Imitations, \emph{Integers} {\bf 9} (2009),
Paper G4, 671-690.
\bibitem[Ray94]{Ray94}J. W. Rayleigh. The Theory of Sound,
\emph{Macmillan, London}, (1894) p. 122-123.
\bibitem[Sko57]{Sko57} Th. Skolem, \"Uber einige Eigenschaften der Zahlenmengen
$[\alpha n+\beta ]$ bei irrationalem $\alpha $ mit einleitenden Bemerkungen
\"uber eine kombinatorishe Probleme, \emph{Norske Vid. Selsk. Forh., Trondheim}
 {\bf 30} (1957), 42-49.
\bibitem[SmSt02]{SmSt02} F. Smith and P. St\u anic\u a, Comply/Constrain
Games or Games with a Muller Twist, \emph{Integers}, {\bf 2}, (2002).
\bibitem[Wyt07]{Wyt07} W.A. Wythoff, A modification of the game of Nim,
\emph{Nieuw Arch. Wisk.} {\bf 7} (1907) 199-202.

\end{thebibliography}
\end{document}